\documentclass[12pt]{article}
\usepackage{txfonts}
\usepackage{mathrsfs}
\usepackage{amssymb}

\begin{document}
\title{A simple proof of the Ohsawa-Takegoshi extension theorem}
\author{Bo-Yong Chen}
\date{}
\maketitle

One of the most beautiful result in complex analysis is the following Ohsawa-Takegoshi extension theorem (cf. [6], see also [1], [3], [5], [7]):

\medskip

\textbf{Theorem.} {\it Let $\Omega$ be a bounded pseudoconvex domain in ${\mathbb C}^n$. Suppose $\sup_\Omega |z_n|^2<e^{-1}$. Then there exists a constant $C_n>0$ such that for every $\varphi\in PSH(\Omega)$,
every holomorphic function $f$ on $\Omega\cap \{z_n=0\}$ with $\int_{\Omega\cap \{z_n=0\}}|f|^2e^{-\varphi}<\infty$, there exists a holomorphic extension $F$ of $f$ to $\Omega$ such that}
$$
\int_\Omega \frac{|F|^2}{|z_n|^2(-\log |z_n|)^2}e^{-\varphi}\le C_n \int_{\Omega\cap \{z_n=0\}}|f|^2e^{-\varphi}.
$$

Recently, there are some attempts to simplify the original proof of Ohsawa-Takegoshi, which is based on a solution of certain {\it twisted} $\bar{\partial}-$equation (cf. [8], [9]). In this paper, we shall give a simple proof by solving directly the $\bar{\partial}-$equation. The idea is inspired by a remarkable paper of Berndtsson-Charpentier (cf. [2]).

\medskip

Let $\rho=\log(|z_n|^2+\epsilon^2)$, $\eta=-\rho+\log(-\rho)$ and $\psi=-\log \eta$, where $\epsilon>0$ is a sufficiently small constant such that $-\rho\ge 1$ on $\Omega$. Since
\begin{equation}
\partial\bar{\partial}\psi=-\frac{\partial\bar{\partial}\eta}\eta+\frac{\partial\eta\bar{\partial}\eta}{\eta^2}=(1+(-\rho)^{-1})\frac{\partial\bar{\partial}\rho}\eta+\frac{\partial\rho\bar{\partial}\rho}{\eta\rho^2}
+\frac{\partial\eta\bar{\partial}\eta}{\eta^2},
\end{equation}
we have $\psi\in PSH(\Omega)$. Put $\phi=\varphi+\log |z_n|^2$. Let $\chi:{\mathbb R}\rightarrow [0,1]$ be a $C^\infty$ cut-off function satisfying $\chi|_{(-\infty,1/2)}=1$ and $\chi|_{(1,\infty)}=0$. By a standard approximation argument we may assume that $f$ is holomorphic in some domain $V\supset\supset\Omega\cap \{z_n=0\}$,  $\varphi$ is $C^\infty$ in a neighborhood of $\overline{\Omega}$, and it suffices to find a holomorphic extension $F$ of $f$ to $\Omega$ such that
$\int_\Omega \frac{|F|^2}{|z_n|^2(-\log |z_n|)^2}e^{-\varphi}\le C_n \int_{V}|f|^2e^{-\varphi}$.
 Thus for $\epsilon$ small enough, we have a well-defined smooth
$\bar{\partial}-$closed $(0,1)$ form $v=f\bar{\partial}\chi(|z_n|^2/\epsilon^2)$ on $\Omega$. Clearly, $v\in L^2_{(0,1)}(\Omega,\phi)$ and there exists a solution of $\bar{\partial}u=v$ with minimal $L^2-$norm in $L^2(\Omega,\phi)$, i.e., $u\bot {\rm Ker\,}\bar{\partial}$. Since $\psi$ is a bounded function, we have $ue^{\psi}\bot {\rm Ker\,}\bar{\partial}$ in $L^2(\Omega,\phi+\psi)$. Thus by H\"ormander's $L^2-$estimates for the $\bar{\partial}$ operator (cf. [4]),
\begin{eqnarray}
\int_\Omega |u|^2 e^{\psi-\phi} & \le & \int_\Omega |\bar{\partial}(ue^\psi)|^2_{\partial\bar{\partial}(\phi+\psi)} e^{-\psi-\phi}\nonumber\\
& = & \int_\Omega |v+\bar{\partial}\psi\wedge u|^2_{\partial\bar{\partial}(\phi+\psi)}e^{\psi-\phi}\nonumber\\
& \le & (1+r^{-1})\int_\Omega |v|^2_{\partial\bar{\partial}(\phi+\psi)}e^{\psi-\phi}+\int_\Omega |\bar{\partial}\psi|^2_{\partial\bar{\partial}(\phi+\psi)}|u|^2e^{\psi-\phi}\nonumber\\
&&+r\int_{{\rm supp\,}v}|\bar{\partial}\psi|^2_{\partial\bar{\partial}(\phi+\psi)}|u|^2e^{\psi-\phi} \ \ \ ({\rm by\ Schwarz's\ inequality})
\end{eqnarray}
where $r>0$ is a small constant to be determined later.  Since $\partial \eta\bar{\partial}\eta=(1+(-\rho)^{-1})^2\partial\rho\bar{\partial}\rho$, we infer from (1) that
$$
\partial\bar{\partial}\psi\ge \frac{\partial\rho\bar{\partial}\rho}{\eta\rho^2}
+\frac{\partial\eta\bar{\partial}\eta}{\eta^2}=\left(\frac1{\eta^2}+\frac1{\eta(-\rho+1)^2}\right)\partial\eta\bar{\partial}\eta.
$$
Thus
\begin{equation}
\int_\Omega |\bar{\partial}\psi|^2_{\partial\bar{\partial}(\phi+\psi)}|u|^2e^{\psi-\phi}\le \int_\Omega \frac{|u|^2}{1+\frac{\eta}{(-\rho+1)^2}} e^{\psi-\phi}.
\end{equation}
By (1), we have $\partial\bar{\partial}\psi\ge \frac{\partial\bar{\partial}\rho}\eta=\frac{\epsilon^2 dz_nd\bar{z}_n}{\eta(|z_n|^2+\epsilon^2)^2}$. Thus by Fubini's theorem, if $\epsilon\ll 1$, we have
\begin{eqnarray}
\int_\Omega |v|^2_{\partial\bar{\partial}(\phi+\psi)}e^{\psi-\phi} &\le & 2\left(\int_{\{\frac{\epsilon^2}2<|z_n|^2<\epsilon^2\}}|\chi'|^2 \frac{(|z_n|^2+\epsilon^2)^2}{\epsilon^2}\frac{|z_n|^2}{\epsilon^4}\frac1{|z_n|^2}\right)\int_{V}|f|^2e^{-\varphi}\nonumber\\
& \le & C_n \int_{V}|f|^2e^{-\varphi}.
\end{eqnarray}
Since
$
\partial\psi\bar{\partial}\psi=\frac1{\eta^2}\left(1+\frac1{-\rho}\right)^2\partial\rho\bar{\partial}\rho\le \frac4{\eta^2}\partial\rho\bar{\partial}\rho
$
and $\partial\bar{\partial}\psi\ge \frac{\epsilon^2 dz_nd\bar{z}_n}{\eta(|z_n|^2+\epsilon^2)^2}\ge \frac{\partial\rho\bar{\partial}\rho}\eta$ on ${\rm supp\,}v$, we get
\begin{equation}
\int_{{\rm supp\,}v} |\bar{\partial}\psi|^2_{\partial\bar{\partial}(\phi+\psi)}|u|^2e^{\psi-\phi}\le \int_\Omega\frac4\eta |u|^2e^{\psi-\phi}.
\end{equation}
Substituting (3),(4),(5) into (2),
$$
\int_\Omega \left(\frac{\frac\eta{(-\rho+1)^2}}{1+\frac\eta{(-\rho+1)^2}}-\frac{4r}\eta\right)|u|^2e^{\psi-\phi}\le (1+r^{-1})C_n \int_{V}|f|^2e^{-\varphi}.
$$
Since $\eta\asymp -\rho$, we may choose  $r=r_n$ sufficiently small such that the left side of the above inequality is bounded below by $c_n\int_\Omega \frac{|u|^2}{|z_n|^2\rho^2}e^{-\varphi}$ for some constant $c_n>0$.  Now put $F_\epsilon=\chi(|z_n|^2/\epsilon^2)f-u$. We conclude that $F_\epsilon$ is a holomorphic extension of $f$ to $\Omega$ together with the estimate
$$
\int_\Omega \frac{|F_\epsilon|^2}{|z_n|^2 (-\log|z_n|^2)^2 }e^{-\varphi}\le C_n'\int_{V}|f|^2e^{-\varphi}.
$$
By taking a weak limit of $F_\epsilon$ as $\epsilon\rightarrow 0$, we get the desired extension.

\medskip

\textbf{Acknowledgement.} The author thanks Dr. Xu Wang for reading the first draft of this note and for his valuable suggestions.

\end{document}